\documentclass[11pt,twoside]{amsart}

\usepackage{mathrsfs}
\renewcommand{\mathcal}{\mathscr}

\usepackage{amsmath,amssymb,amstext}
\usepackage{hyperref}
\newtheorem{theorem}{Theorem}

\newtheorem{lemma}[theorem]{Lemma}

\newtheorem{preremark}[theorem]{Remark}

\newcommand{\N}{{\mathbb N}}
\newcommand{\R}{{\mathbb R}}

\newcommand{\AAA}{{\mathcal{A}}}

\newcommand{\VVV}{{\mathcal{V}}}

\newcommand{\e}{{\varepsilon}}

\renewcommand{\geq}{\geqslant}
\renewcommand{\ge}{\geqslant}
\renewcommand{\leq}{\leqslant}
\renewcommand{\le}{\leqslant}

\begin{document}

\author{Yannick Sire}
\address{{\em YS} --
Universit\'e Aix-Marseille 3, Paul C\'ezanne --
LATP --
Marseille, France}
\email{sire@cmi.univ-mrs.fr}

\author{Enrico Valdinoci}
\address{{\em EV} --
Universit\`a di Roma Tor Vergata --
Dipartimento di Matematica --
I-00133 Rome, Italy}
\email{enrico@mat.uniroma3.it}

\thanks{EV has been  
supported
by
FIRB
project ``Analysis and Beyond''
and GNAMPA
project
``Equa\-zio\-ni non\-li\-nea\-ri su va\-rie\-t\`a:
pro\-prie\-t\`a qua\-li\-ta\-tive e clas\-si\-fi\-ca\-zio\-ne
del\-le so\-lu\-zio\-ni''. 
Part of this work
was carried out
while YS was visiting Tor Vergata.}

\title{Density estimates for phase transitions with a trace}

\begin{abstract}
We consider a functional obtained by adding a trace term to
the Allen-Cahn phase segregation
model and we prove some density estimates for
the level sets of the interfaces.

We treat in a unified way also the cases of
possible degeneracy and singularity of the ellipticity
of the model and the quasiminimal case.
\end{abstract}

\maketitle

\section{Introduction}

Let~$p\in(1,+\infty)$,
$n\ge2$ and $\Omega$ be an open bounded subset 
of~$\R^n$.
For any~$u\in W^{1,p}_{\rm loc}(\R^n,[-1,1])$,
we define the functional
\begin{equation}\label{F}
{\mathcal{E}}_\Omega (u):=\int_{\Omega\cap\R^n_+} |\nabla
u(x)|^p+F(u(x))\,dx
+\int_{ \Omega\cap \{ x_n=0\}} G(u(x',0))\,dx',\end{equation}
where we used the notation~$x:=(x',x_n)\in \R^{n-1}\times\R$
and $\R^n_+:=\R^{n-1}\times(0,+\infty)$.
Also, for any~$R>0$ and any~$x\in \R^n$, we denote by
$B_R^n(x)$ the Euclidean,
open, $n$-dimensional ball centered at~$x$, 
and~$B_R^n:=B^n_R(0)$.
We set~$B_R^+(x):=B_R^n(x)\cap
\R^n_+$, $B^+_R:=B^+_R(0)$, and
we use the short notation
$$ {\mathcal{E}}_{R,x_o} (u):=
{\mathcal{E}}_{B_R(x_o)}(u)=
\int_{B_R^+(x_o)} |\nabla 
u(x)|^p+F(u(x))\,dx
+\int_{ B_R^{n-1}(x_o)} G(u(x',0))\,dx'.$$
We will suppose that~$F$
and~$G$ are non-negative ``double-well'' potentials.
More precisely, and in fact more
generally,
we assume that there exists~$C_o\ge 1$ 
such 
that,
for any~$\tau\in\R$, we have
\begin{equation}\label{d9832ooo2}
\max\{ F(\tau), G(\tau)\}\le{C_o} (1-\tau^2)^p {\mbox{ and }}
F(\tau)\ge \frac{1}{C_o} (1-\tau^2)^p.\end{equation}
A paradigmatic example is given by $F(\tau)=G(\tau)=(1-\tau^2)^p$,
but more general potentials are allowed by \eqref{d9832ooo2}.
The gradient term in~\eqref{F} is reminiscent of a $p$-Laplacian
partial differential equation (hence, it
encodes a possibly singular or degenerate ellipticity).
We remark that the functional in~\eqref{F} reduces to
the standard Allen-Cahn phase segregation
model when~$G$ is identically zero and~$\Omega$ lies
in~$\{x_n>0\}$. Thus, in a sense, the functional in~\eqref{F}
represents a phase transition in~$\R^n_+$
with a double-well~$G$ keeping track of a phase segregation
on the trace of~$\Omega$ along~$\{x_n=0\}$
and it may be seen as a toy-model to understand
the more complicated phenomena
arising in non-local phase transitions, which
have been the object of an extensive study in recent years
(see, among the others,~\cite{ABS94, AB98, ABS98, Mar}
and also~\cite{Sol, SV, SV2, Cinti} for a relation between fractional
operators and boundary reactions).
In practical situations,
the non-local effects may be the consequence
of a long-range
interaction between particles, as it happens
in some statistical mechanics models
(see, e.g.,~\cite{Orl}).

The trace term $\int_{ \Omega\cap \{ x_n=0\}} G(u(x',0))\,dx'$
may also be considered as a model for taking
into account the effect of the boundary of the
container in which the phase transition occurs:
in this framework, the container is $\R^n_+$,
which, of course, up to a blow up, is a simplified,
but effective,
version of a smooth container
when we are interested in the behavior near
its boundary.
In this sense, we hope
that this paper
may be as a first step towards a
more comprehensive study of the geometric features
of the phase transitions under
even more severe boundary
and non-local effects.

Given~$Q\ge1$, we say that $u$ is a $Q$-minimizer 
if~${\mathcal{E}}_\Omega(u)<+\infty$
and
\begin{equation}\label{QM} \begin{split}&{\mathcal{E}}_\Omega(u)\le
Q {\mathcal{E}}_\Omega(u+\varphi)\\ &\qquad{\mbox{
for all bounded and open~$\Omega
\subset\R^n$}}
\\ &\qquad{\mbox{
and all Lipschitz continuous
functions~$\varphi$ supported in $\Omega$.}}\end{split}
\end{equation}
The case of $Q$-minimizers in a fixed domain $\Omega_o$
may be treated in a similar way (just suppose that $\Omega\subseteq
\Omega_o$ in \eqref{QM} and so on).
The study of~$Q$-minimizers is a classical topic in the
calculus of variations (see, e.g.,~\cite{GG}).
When~$Q=1$ in \eqref{QM}, $u$ is usually said to be a minimizer.
It is easily seen that when $p=2$
the minimizers satisfy the partial differential equation
problem with Neumann condition 
$$ \left\{\begin{matrix} 2\Delta u= F'(u) & {\mbox{in $\R^n_+$,}}\\
2\partial_{x_n} u=G'(u)
&{\mbox{ on $\{x_n=0\}$.}}
\end{matrix}
\right.$$
Such type of problems have been studied in \cite{Sol, SV}.
Analogously, the minimizers for $p\in (1,2)\cup(2,+\infty)$
satisfy a quasilinear partial differential equation whose
ellipticity becomes singular or degenerate at the critical points
of the solution, and the corresponding Neumann condition
becomes non-linear too: these types of problem have been studied,
for instance, in \cite{SV2}.\medskip

This is the main result of this paper:

\begin{theorem}\label{1}
Let~$\mathcal{L}^n$
denote the~$n$-dimensional Lebesgue measure.
Let~$u$ be a Lipschitz continuous~$Q$-minimizer.

Then, there exists a positive~$C_*$, only
depending on~$n$, $Q$, $p$,
the quantity~$C_o$ in~\eqref{d9832ooo2}
and the Lipschitz constant of~$u$,
such that
\begin{equation}\label{DS}
{\mathcal{E}}_{R,x_o}(u)\le C_* R^{n-1}\end{equation}
for any~$x_o\in \overline{\R^n_+}$ and any $R\ge 1$.

Furthermore, given any~$\theta\in(-1,1)$, if we suppose that
there exist two positive real numbers~$\mu_1$
and~$\mu_2$ such that
\begin{equation}
\label{eqpass0}
\mathcal{L}^n\Big(B_{\mu_1}^+(x_o)\cap\{ u>\theta\}\Big)
\geq\mu_2,
\end{equation} 
then
there exist positive~$r_0$
and~$c$, which depend 
only on~$n$, $Q$, $p$, $\theta$,~$\mu_1$,~$\mu_2$, 
the quantity~$C_o$ in~\eqref{d9832ooo2}
and the Lipschitz constant of~$u$, in
such a way that
\begin{equation}\label{DS1}
\mathcal{L}^n\Big(B_r^+(x_o)\cap\{ u>\theta\}\Big) \geq c\,r^n,
\end{equation} 
for
any~$r\geq r_0$.

Analogously, if
\begin{equation}\label{eqpass0-II}
\mathcal{L}^n\Big(B_{\mu_1}^+(x_o)\cap\{ 
u<\theta\}\Big)
\geq\mu_2\end{equation}
then
\begin{equation}\label{DS2}
\mathcal{L}^n\Big(B_r^+(x_o)\cap\{ u<\theta\}\Big) \geq c\,r^n,
\end{equation}
for
any~$r\geq r_0$.
\end{theorem}

We remark that \eqref{eqpass0} (respectively,
\eqref{eqpass0-II}) is satisfied if $u(x_o)>\theta$ (respectively,
$u(x_o)<\theta$): in this case, $\mu_1$ and $\mu_2$
just depend on $|\theta-u(x_o)|$ and on the modulus
of continuity of $u$.

Our Theorem~\ref{1} fits into the line
of research of density estimates for
phase transition, as started in~\cite{CC}, to which
it reduces when $G:=0$ or, basically, when we
look at balls $B_r(x_o)$ that do not
intersect $\{ x_n=0\}$. Namely,
the purpose of this type of
researches is to try to understand
how the level sets of a ``good'' solution~$u$
behave in measure. Such level sets are physically
very relevant, since they represent, roughly
speaking, the separation interface of the two phases~$+1$ and~$-1$
in the Allen-Cahn system.
Also, from these measure theoretic estimates, it
is possible to deduce a locally uniform convergence of the level 
sets at the~$\Gamma$-limit, and this information
plays a crucial role in some rigidity problems
(see~\cite{S-T, VSS, S-A}).

Among the many extensions of~\cite{CC}, we recall here
the ones in~\cite{PV}, where the~$p$-Laplacian case has
been considered, \cite{FV}, for quasiminima,
and the recent preprint~\cite{OSV}, dealing with
a fully non-local case.

Estimate \eqref{DS1}
(as well as \eqref{DS2}) is obviously optimal (up to the
constant $c$), because of the trivial upper bound
$$ \mathcal{L}^n\Big(B_r^+(x_o)\cap\{ u>\theta\}\Big)\leq
\mathcal{L}^n\Big(B_r^+(x_o)\Big)\le
\mathcal{L}^n\Big(B_r(x_o)\Big) \sim r^n.$$
Estimate \eqref{DS} is optimal too,
as shown by the case $G:=0$, taking as $u(x)=u_o(\omega\cdot x)$,
where $\omega\in {\rm S}^{n-1}$ and $u_o:\R\rightarrow\R$ is a
minimizer of the one-dimensional Allen-Cahn functional.

As far as we know, in the framework of
the functional in~\eqref{F}, Theorem~\ref{1} of this
paper is new
even in the cases~$p=2$ (i.e., when
the diffusion term reduces to the standard Laplacian) and~$Q=1$
(i.e., for minimizers).

\section{Proof of Theorem \ref{1}}

We will denote by~``${\,\rm const\,}$'' suitable positive
constants (possibly different line by line) only depending
on the quantities fixed in the hypotheses of Theorem \ref{1}.
First of all, we prove~\eqref{DS}. This is done
by a technique well-developed after~\cite{CC}: given
any~$x_o\in \overline{\R^n_+}$ and any $R\ge 1$, we
take~$\beta\in
C^\infty(\R^n)$, with~$\beta(x)=-1$ for any~$x\in B_{R-(1/2)}(x_o)$
and~$\beta(x)=1$ for any~$x\in \R^n-B_{R-(1/4)}(x_o)$,
with~$|\nabla\beta(x)|\le 50$ for any~$x\in\R^n$.
Let~$w(x):=\min\{u(x),\beta(x)\}$. Then, $w(x)=u(x)$
in~$\R^n-B_R(x_o)$, 
and $w=-1$ in $B_{R-(1/2)}(x_o)$.
So, by~\eqref{QM},
we obtain that
\begin{equation}\label{791}
\begin{split}
& \frac1Q
{\mathcal{E}}_{R,x_o} (u)\le{\mathcal{E}}_{B_R(x_o)}(w)\\ 
&\qquad=
\int_{(B_R^n (x_o)- B^n_{R-(1/2)}(x_o))\cap\R^n_+ } |\nabla
w(x)|^p+F(w(x))\,dx
\\ &\qquad\qquad+\int_{ (B_R^{n}(x_o)- B^{n}_{R-(1/2)}(x_o))\cap\{x_n=0\}
} 
G(w(x',0))\,dx'.
\end{split}\end{equation}
Moreover, we have that
$|\nabla w(x)|\le
|\nabla u(x)|+|\nabla\beta(x)|\le{\,\rm const\,}$, and so~\eqref{791}
gives that
\begin{eqnarray*}&& {\mathcal{E}}_{R,x_o} (u) 
\le {\,\rm const\,} \left[ 
{\mathcal{L}}^n
\Big( {(B_R^n(x_o) - B^n_{R-(1/2)}(x_o))\cap\R^n_+ }\Big)\right.
\\ &&\qquad\quad\left.+
{\mathcal{H}}^{n-1} \Big(
(B_R^{n}(x_o)- B^{n}_{R-(1/2)}(x_o))\cap\{x_n=0\} 
\Big)\right]\\ &&\qquad\le 
{\,\rm const\,} R^{n-1},\end{eqnarray*}
where we denoted by~${\mathcal{H}}^{n-1}$ the~$(n-1)$-dimensional
Hausdorff measure. This
proves~\eqref{DS}.\medskip

Now, we prove~\eqref{DS1} (the proof of~\eqref{DS2}
is the same and it will be omitted).
The proof of~\eqref{DS1} that
we give here
is a modification of one of the
proofs performed in~\cite{FV}, which was inspired 
by~\cite{S}
(other approaches, as the ones in \cite{CC, PV}
are also possible, but they may require
additional assumptions). Differently from the
existing literature, here
some technical complications arise
in order to cope with the trace term of the functional along~$\{x_n=0\}$.
Indeed, even if such a term behaves as an~$(n-1)$-dimensional
correction, and therefore
may look negligible,
it is not completely clear that it does not dangerously interact
with some ``area terms'' arising in the density estimates,
such as the 
bound in~\eqref{DS} and the subsequent quantities in~\eqref{FVp12}.
For this, we will have to perform some careful computation.

First, we observe that once~\eqref{DS1} is proved for
some~$\theta_o\in (-1,-1/2]$, then it is proved for 
all~$\theta\in[\theta_o,1)$, because
$${\mathcal{E}}_{r,x_o}(u)\ge\int_{B_r^+(x_o)}F(u)\,dx
\ge \inf_{[\theta_o,\theta]} F\,\cdot\,
\mathcal{L}^n\Big(B_r^+(x_o)\cap\{ \theta>u>\theta_o\}\Big),$$
and
if~\eqref{eqpass0} holds for~$\theta\in[\theta_o,1)$, it holds
for~$\theta_o$ too, so using~\eqref{DS1} for~$\theta_o$
and~\eqref{DS} we obtain
\begin{eqnarray*}
&& \mathcal{L}^n\Big(B_r^+(x_o)\cap\{ u>\theta\}\Big) 
\ge
\mathcal{L}^n\Big(B_r^+(x_o)\cap\{ u>\theta_o\}\Big)  
-
\mathcal{L}^n\Big(B_r^+(x_o)\cap\{ \theta>u>\theta_o\}\Big)\\
&&\qquad\ge
cr^n -\,{\rm const\,} {\mathcal{E}}_{r,x_o} (u) \ge
cr^n -\,{\rm const\,} r^{n-1} \ge \frac{c}{2}r^n
\end{eqnarray*}
if~$r\ge r_o$ and~$r_o$ is large enough (here the
``${\,\rm const\,}$'' may depend on the fixed~$\theta_o$ too). This would 
be the
proof of~\eqref{DS1} for any~$\theta\in[\theta_o,1)$, up
to relabeling~$c$, and therefore, in what follows,
we will assume, without any restriction,
that
\begin{equation}\label{RES}
\theta\in(-1,-1/2].\end{equation}
Moreover,
we observe that the portion of space~$\R^n\cap \{x_n<0\}$
does not play any role in Theorem~\ref{1}, in the sense that,
if we define
$$ \tilde u(x)=\tilde u(x',x_n):=\left\{\begin{matrix}
u(x',x_n) & {\mbox{ if $x_n\ge0$,}}\\
u(x',-x_n) & {\mbox{ if $x_n<0$,}}
\end{matrix}\right. $$
we have that~$\tilde u$ is Lipschitz, since 
so is $u$, that ${\mathcal{E}}_{\Omega}(\tilde 
u+\varphi)=
{\mathcal{E}}_{\Omega}(u+\varphi)$ for any 
perturbation~$\varphi$ in \eqref{QM},
that~$\tilde u$ is a $Q$-minimizer and that if~\eqref{DS1}
holds for~$\tilde u$ then it holds for~$u$ as well.
Consequently, we replace~$u$ with~$\tilde u$ and then we drop
the superscript tilde, that is we may and do suppose that
\begin{equation}\label{symmetry}
u(x',-x_n)=u(x',x_n).
\end{equation}
This symmetry property will play an important role,
by allowing us to disregard some trace term
in a subsequent isoperimetric inequality (that is
\eqref{ISO} below: roughly speaking, this trick will
make the trace term
in the density estimates always be weighted by
the potential, thus killing any unweighted geometric measure
on~$\{x_n=0\}$).

We 
take
\begin{equation}\label{ET}
\begin{split}
&{\mbox{$T$ to be a
free parameter, that in the sequel will be chosen
to be 
suitably large,}}\\&\qquad{\mbox{
possibly in dependence of the quantities fixed
in the statement of Theorem~\ref{1},}}
\\&\qquad{\mbox{
and also in dependence of a further auxiliary 
parameter $\varepsilon$}}
\\&\qquad{\mbox{
that will
be introduced later on, after \eqref{6dwrdrrq00wk1}.}}
\end{split}\end{equation}
We
set
$$ S(\tau):=\min \left\{
{(\tau+1)^p},\,1
\right\} {\mbox{ for any
$\tau\in\R$.}}$$
Also, for any $x\in \overline{\R^n_+}$
and any $k\in\N$, we let
$$ v_k (x):= \left\{ \begin{matrix} 2 e^{|x-x_o|-(k+1)T}-1 &
{\mbox{ for any
$x\in B_{(k+2)T}(x_o)\cap \overline{\R^n_+}$,}}\\
2e^T-1& {\mbox{ for any
$x\in \overline{\R^n_+}-B_{(k+2)T}(x_o)$.}}\end{matrix}\right.
$$
If $x\in \R^n\cap \{x_n<0\}$, we also define
\begin{equation}\label{symmetry2}
v_k(x',x_n):=v_k(x',-x_n).
\end{equation}
By construction,
$v_k$ is Lipschitz.
Furthermore, we deduce from~\eqref{d9832ooo2}
that
\begin{equation}\label{Sa0}
\begin{split}
|\nabla v_k(x)|^p =& (2 e^{|x-x_o|-(k+1)T})^p\\
=& \,(v_k(x)+1)^p\\
\leq& {\,\rm const\,} S(v_k(x))
\,\end{split}
\end{equation}
for any $x\in B^+_{(k+2)T}(x_o)$, and therefore,
by \eqref{symmetry2}, for almost every $x\in\R^n$.
Furthermore, we see
that
\begin{equation}\label{62bis} \max\{F(\tau), G(\tau)\} \leq
{\,\rm const}\, S(\tau)\end{equation}
for any~$\tau\in[-1,1]$, and that
\begin{equation}\label{61}
F(\tau)\geq \,{\rm const}\,(\tau+1)^p
=\,{\rm const}\,S(\tau)
\end{equation}
when~$\tau\in [-1,-1/2]$.

We remark that
\begin{equation}\label{VK}
{\mbox{if $x\in\R^n_+$ and
$|x-x_o|>(k+1)T$, then $v_k(x)>1\ge u(x)$
}}\end{equation}
and so, recalling \eqref{symmetry}
and \eqref{symmetry2},
we conclude that~$\{u>v_k\}=
\{ x\in\R^n {\mbox{ s.t. }} u(x)>v_k(x) \}$
is a bounded set.
Accordingly, we can make use of~\eqref{QM} with~$\Omega:=\{u>v_k\}$.
This, and the use of~\eqref{Sa0}
and~\eqref{62bis},
imply the following estimate:
\begin{equation}\label{Sa1}
\begin{split}
&\int_{\{ u>v_k \}\cap\R^n_+} |\nabla u|^p+F(u)\,dx+
\int_{\{ u>v_k \}\cap\{x_n=0\} } G(u)\,dx'
\\&\qquad ={\mathcal{E}}_{ \{ u>v_k \} }(u)
\\&\qquad \le Q {\mathcal{E}}_{ \{ u>v_k \} }(v_k)
\\&\qquad =Q\left[
\int_{\{ u>v_k \}\cap\R^n_+} |\nabla v_k|^p+F(v_k)\,dx
+\int_{\{ u>v_k \}\cap\{x_n=0\} } G(v_k)\,dx'
\right]
\\&\qquad \leq {\,\rm const\,}\left[
\int_{\{ u>v_k \} \cap\R^n_+} S(v_k)\,dx+
\int_{\{ u>v_k \}\cap\{x_n=0\} } S(v_k)\,dx'
\right]
.\end{split}
\end{equation}
Now, we make a general observation:
given any Lipschitz function~$w$
on a measurable set~$U\subseteq \R^n$ with image in~$[-1,1]$,
we have 
\begin{equation}\label{fsgh6710001}
\begin{split}
\int_U |\nabla w|^p+F(w)\,dx
\geq&{\;\rm const\,}
\int_U |\nabla w| \Big(
F(w)\Big)^{(p-1)/p}
\,dx
\\=&{\,\rm const\,}
\int_{-1}^1
\Big(F(\tau)\Big)^{(p-1)/p}
{\mathcal{H}}^{n-1}\Big( U\cap
\{ w=\tau\}\Big)
\,d\tau
,\end{split}
\end{equation}
due to
the Young inequality and the coarea formula.

Also, we define
\begin{eqnarray*}
{\mathcal{M}}_k(\tau)&:=&\{x\in\R^n {\mbox{ s.t. }}
\tau=u(x)\ge v_k(x)
\}\\ &=&
\{u\geq v_k\}\cap\{u=\tau\}
\\{\mbox{and }}\quad
{\mathcal{N}}_k(\tau)&:=&\{x\in\R^n {\mbox{ s.t. }}   
u(x)\ge v_k(x)=\tau
\}\\ &=&
\{u\geq v_k\}\cap
\{ v_k=\tau\},
\end{eqnarray*}
and we remark that
\begin{eqnarray*} && B_{kT}^+(x_o)\cap\{ u>\theta\}\subseteq
\{ x\in B_{kT}^+(x_o) {\mbox { s.t. } }
u(x)>\tau>v_k(x)\}
\\ &&\subseteq
\{ x\in \R^n {\mbox { s.t. } }
u(x)>\tau>v_k(x)\} = \{u>\tau>v_k\}
\end{eqnarray*}
for any~$\tau\in [(\theta-1)/2,\theta]$
as long as the free parameter~$T$
is chosen large enough. 

We employ the latter formula and
the isoperimetric inequality to obtain
\begin{equation}\label{ISO}
\begin{split}
& \Big({\mathcal{L}}^n (
B_{kT}(x_o)\cap\{ u>\theta\})\Big)^{(n-1)/n}\leq
\Big({\mathcal{L}}^n (
\{u>\tau>v_k\}
)\Big)^{(n-1)/n}\\ \qquad
&\leq
{\rm const\,} \Big(
{\mathcal{H}}^{n-1}({\mathcal{M}}_k(\tau))+
{\mathcal{H}}^{n-1}({\mathcal{N}}_k(\tau))\Big)
,\end{split}
\end{equation}
for any~$\tau\in [(\theta-1)/2,\theta]$.

Notice that we have in the back of our mind here the symmetry
in~\eqref{symmetry} and~\eqref{symmetry2},
since we are willing to estimate sets in~\eqref{ISO}
in the whole
of~$\R^n$ instead of~$\R^n_+$: due to such a symmetry,
this choice will be paid only by a factor~$2$ later
on: see~\eqref{90}. On the contrary,
without this trick we would have got
also a term of the form~${\mathcal{H}}^{n-1}(B^{n-1}_{kT})$
in~\eqref{ISO}, and this
would have risked to be too large to be controlled
by the quantities in~\eqref{DS}
and~\eqref{FVp12}.

Making use of~\eqref{ISO} and then
of~\eqref{fsgh6710001} with~$U:=\{u\ge 
v_k\}$
and either~$w:=u$ or~$w:=v_k$,
we conclude that
\begin{equation}\label{6dwrdrrq00wk0-0}
\begin{split}
&\Big({\mathcal{L}}^n (
B_{kT}^+(x_o)\cap\{ u>\theta\})\Big)^{(n-1)/n}\\
&= {\rm const\,} \int_{(\theta-1)/2}^{\theta}
\Big(F(\tau)\Big)^{(p-1)/p}\,d\tau
\,\Big({\mathcal{L}}^n (
B_{kT}(x_o)\cap\{ u>\theta\})\Big)^{(n-1)/n}
\\&\leq
{\rm const\,} \int_{-1}^{1}
\Big( F(\tau)\Big)^{(p-1)/p}
\Big(
{\mathcal{H}}^{n-1}\big({\mathcal{M}}_k(\tau)\big)
+{\mathcal{H}}^{n-1}\big({\mathcal{N}}_k(\tau)\big)\Big)
\,d\tau
\\ &\leq
{\rm const\,} \left(
\int_{
\{u\ge v_k\}
} |\nabla u|^p+ F(u)\,dx
+
\int_{
\{u\ge
v_k\}
} |\nabla v_k|^p+F(v_k)\,dx
\right).\end{split}\end{equation}
Now, we remark that
\begin{equation}\label{90}
\begin{split}
&\int_{
\{u\ge v_k\}
} |\nabla u|^p+ F(u)\,dx
+
\int_{
\{u\ge
v_k\} 
} |\nabla v_k|^p+F(v_k)\,dx
\\ &\qquad = 2\left[
\int_{
\{u\ge v_k\}\cap\R^n_+
} |\nabla u|^p+ F(u)\,dx
+
\int_{
\{u\ge
v_k\} \cap\R^n_+
} |\nabla v_k|^p+F(v_k)\,dx
\right],
\end{split}
\end{equation}
thanks to~\eqref{symmetry} and~\eqref{symmetry2}.

Therefore, exploiting~\eqref{6dwrdrrq00wk0-0},
\eqref{90},
\eqref{Sa0},
\eqref{62bis}  
and~\eqref{Sa1},
we obtain
\begin{equation*}
\begin{split}
&\Big({\mathcal{L}}^n (
B_{kT}^+(x_o)\cap\{ u>\theta\})\Big)^{(n-1)/n}\\
&\leq {\rm const\,}
\left[
\int_{
\{u\ge v_k\}\cap\R^n_+
} |\nabla u|^p+ F(u)\,dx
+
\int_{
\{u\ge
v_k\} \cap\R^n_+
} |\nabla v_k|^p+F(v_k)\,dx
\right]
\\ &\leq
{\rm const\,}\left[
\int_{\{ u\ge
v_k \}\cap\R^n_+} S(v_k)\,dx
+
\int_{\{ u\ge
v_k \}\cap\{x_n=0\}} S(v_k)\,dx'
\right].
\end{split}
\end{equation*}
That is, recalling~\eqref{VK},
\begin{equation}\label{6dwrdrrq00wk0}
\begin{split}
&\Big({\mathcal{L}}^n (
B_{kT}^+(x_o)\cap\{ u>\theta\})\Big)^{(n-1)/n}\\
&\leq
{\rm const\,}\left[
\int_{\{ u\ge
v_k \}\cap B^+_{(k+1)T} (x_o)} S(v_k)\,dx
+
\int_{\{ u\ge
v_k \}\cap
B_{(k+1)T}(x_o)\cap
\{x_n=0\}} S(v_k)\,dx'
\right].
\end{split}
\end{equation}
We define
\begin{eqnarray*}
&& \ell_1=\ell_1(k):=\int_{B_{kT}^+(x_o)} S(v_k)\,dx
\\ &&
\ell_2=\ell_2(k):=\int_{ \{u\ge v_k
\}\cap (B_{(k+1)T}^+(x_o)-B_{kT}^+(x_o))} 
S(v_k)\,dx
\\ \,{\mbox{ and }}\,&&
\ell_3=\ell_3(k):=
\int_{\{ u\ge
v_k \}\cap B_{(k+1)T}(x_o)\cap\{x_n=0\}} S(v_k)\,dx'.
\end{eqnarray*}
With this notation,
we see that~\eqref{6dwrdrrq00wk0}
can be written as
\begin{equation}\label{6dwrdrrq00wk1}
\Big({\mathcal{L}}^n (
B_{kT}^+(x_o)\cap\{ u>\theta\})\Big)^{(n-1)/n}\leq\,{\rm const}\,
(\ell_1+\ell_2+\ell_3).
\end{equation}
Now, we fix a small~$\varepsilon>0$
to be taken appropriately small (in fact, at the end, this~$\varepsilon$
will be fixed explicitly in Lemma~\ref{H12}
below) and we claim that
\begin{equation}\label{L3}
\ell_3 \le \,{\rm const}\,
\varepsilon k^{n-1}+
C_\varepsilon (kT)^{n-2},\end{equation}
for a suitable $C_\varepsilon>0$.
The proof of \eqref{L3} is indeed a bit long and complicated
and it will be completed only below \eqref{775}, after
some delicate computations.
To prove \eqref{L3}, first we notice that
when $|x_{o,n}|>
(k+1)T$ then
$B_{(k+1)T}(x_o)\cap\{x_n=0\}=\varnothing$, so $\ell_3=0$
and \eqref{L3} is obviously true. As a consequence, we may suppose
that 
$$ |x_{o,n}|\le
(k+1)T$$
and so
we can define \begin{equation}\label{R5}
\rho_k:= \sqrt{(k+1)^2 T^2 
-x_{o,n}^2}.\end{equation}
We see that
$$ B^n_{(k+1)T}(x_o) \cap \{ x_n=0\}
\subseteq B^{n-1}_{\rho_k}(x'_o)$$
and therefore
\begin{equation}\label{L3-cor-ii}
\begin{split}
& \ell_3 \le{\,\rm const\,}\int_{B^{n}_{(k+1)T}(x_o) \cap \{
x_n=0\}}
(v_k(x',0)+1)^p\,dx'
\\ &\qquad \le {\,\rm const\,}
\int_{B^{n-1}_{\rho_k}(x'_o)} e^{
p\big( \sqrt{|x'-x'_o|^2+x_{o,n}^2}-(k+1)T\big)
} \,dx'\\&\qquad = {\,\rm const\,} e^{ -p(k+1)T }
\int_0^{\rho_k} r^{n-2}
e^{
p \sqrt{r^2+x_{o,n}^2}
} \,dr.
\end{split}\end{equation}
Now, to prove \eqref{L3},
we distinguish two cases: either $n\ge3$
or $n=2$.

If $n\ge3$, we make use of \eqref{L3-cor-ii}
to conclude that
\begin{equation}\label{L3-cor-iii}
\begin{split}
&\ell_3\le {\,\rm const\,} \rho_k^{n-3} e^{ -p(k+1)T }
\int_0^{\rho_k} r
e^{
p \sqrt{r^2+x_{o,n}^2}
} \,dr,\end{split}\end{equation}
and we perform the substitution
\begin{equation}\label{S34}
s:= \sqrt{r^2+x_{o,n}^2}.\end{equation}
We obtain that
$ s\,ds = r \,dr$, hence
\eqref{L3-cor-iii}
becomes
\begin{equation*}
\begin{split}
& \ell_3 \le
{\,\rm const\,} \rho_k^{n-3} e^{ -p(k+1)T }
\int_{|x_{o,n}|}^{(k+1)T} s
e^{
p s
} \,ds
\\ &\qquad \le {\,\rm const\,}
\rho_k^{n-3} \, (k+1)\,T\, e^{ -p(k+1)T }
\int_{|x_{o,n}|}^{(k+1)T}
e^{
p s
} \,ds
\\ &\qquad \le {\,\rm const\,}
\rho_k^{n-3} \, (k+1)\,T
\\ &\qquad \le {\,\rm const\,}
\big( (k+1) T\big)^{n-2}
\\ &\qquad \le {\,\rm const\,}
(kT)^{n-2}.\end{split}\end{equation*}
This proves \eqref{L3} when $n\ge3$, so
now we deal with the proof of \eqref{L3}
when $n=2$: in this case, we claim that
\begin{equation}\label{L3ux}
\int_0^{\rho_k}e^{
p \sqrt{r^2+x_{o,n}^2}   
} \,dr \le \big( 2\varepsilon^2 kT+\tilde C_\varepsilon\big)
e^{p(k+1)T}
\end{equation}
for a suitable $\tilde C_\varepsilon>0$.

To prove \eqref{L3ux},
we distinguish two sub-cases:
either $\rho_k<\varepsilon^2 (k+1)T$
or $\rho_k\ge\varepsilon^2 (k+1)T$.

If $\rho_k<\varepsilon^2 (k+1)T$, we have that
\begin{eqnarray*}
&& \int_0^{\rho_k}
e^{
p \sqrt{r^2+x_{o,n}^2}
} \,dr \le e^{
p \sqrt{\rho_k^2+x_{o,n}^2}
} \rho_k = e^{p(k+1)T}\rho_k
\\ &&\qquad\le e^{p(k+1)T}
\varepsilon^2 (k+1)T\le  2e^{p(k+1)T}      
\varepsilon^2 kT,
\end{eqnarray*}
and this proves \eqref{L3ux} in the sub-case
$\rho_k<\varepsilon^2 (k+1)T$.

Now, we prove \eqref{L3ux} in the sub-case
$\rho_k\ge\varepsilon^2 (k+1)T$, that gives, recalling \eqref{R5},
\begin{equation}\label{3Y}
|x_{o,n}|\le \sqrt{1-\varepsilon^4}\,(k+1)T.
\end{equation}
Then, we make
the substitution in \eqref{S34},
we split the domains of integration, and we obtain
\begin{equation}\label{XY}
\begin{split}
&\int_0^{\rho_k} 
e^{
p \sqrt{r^2+x_{o,n}^2}   
} \,dr =
\int_{|x_{o,n}|}^{(k+1)T}   
e^{ps} \frac{s}{\sqrt{s^2 -x_{o,n}^2}}\,ds\\
&\qquad\le
\int_{|x_{o,n}|}^{{\sqrt{1+\varepsilon^4} {|x_{o,n}|} }} 
e^{ps} \frac{s}{\sqrt{s^2 -x_{o,n}^2}}\,ds
+\Xi
\int_{
{\sqrt{1+\varepsilon^4} {|x_{o,n}|} }
}^{(k+1)T} 
e^{ps} \frac{s}{\sqrt{s^2 -x_{o,n}^2}}\,ds,
\end{split}
\end{equation}
where
$$ \Xi:=\left\{
\begin{matrix}
1 & {\mbox{if $
{{\sqrt{1+\varepsilon^4} {|x_{o,n}|} }
} < {(k+1)T} $,}}\\
0 & {\mbox{if $    
{{\sqrt{1+\varepsilon^4} {|x_{o,n}|} }
} \ge {(k+1)T} $.}}
\end{matrix}
\right.$$
So, we compute separately the latter two integrals
in \eqref{XY}.
For the first one, 
recalling 
\eqref{3Y}, we have:
\begin{equation}\label{XY1}\begin{split}
& \int_{|x_{o,n}|}^{{\sqrt{1+\varepsilon^4} {|x_{o,n}|} }}  
e^{ps} \frac{s}{\sqrt{s^2 -x_{o,n}^2}}\,ds
\le
e^{p{\sqrt{1+\varepsilon^4} {|x_{o,n}|} }}
\int_{|x_{o,n}|}^{{\sqrt{1+\varepsilon^4} {|x_{o,n}|} }}
\frac{s}{\sqrt{s^2 -x_{o,n}^2}}\,ds
\\ &\qquad
=
e^{p{\sqrt{1+\varepsilon^4} {|x_{o,n}|} }}
\varepsilon^2 {|x_{o,n}|}
\le
e^{p{\sqrt{1-\varepsilon^8} {(k+1)T} }}
\varepsilon^2 (k+1)T\le 2e^{p{{(k+1)T} }}   
\varepsilon^2 kT.
\end{split}\end{equation}
Now we estimate
the last integral in \eqref{XY} as follows:
\begin{equation}\label{XY2}
\begin{split}
& \Xi
\int_{
{\sqrt{1+\varepsilon^4} {|x_{o,n}|} }
}^{(k+1)T}
e^{ps} \frac{s}{\sqrt{s^2 -x_{o,n}^2}}\,ds
=
\Xi
\int_{
{\sqrt{1+\varepsilon^4} {|x_{o,n}|} }
}^{(k+1)T}
e^{ps} \sqrt{ 1+\frac{x_{o,n}^2}{{s^2 -x_{o,n}^2}} }\,ds
\\ &\qquad
\le 
\Xi
\int_{
{\sqrt{1+\varepsilon^4} {|x_{o,n}|} }
}^{(k+1)T}
e^{ps} \sqrt{ 1+\frac{1}{{ \varepsilon^4}} }\,ds     
\le \sqrt{ 1+\frac{1}{{ \varepsilon^4}} }
\int_{-\infty}^{(k+1)T}   
e^{ps} \,ds   
=\tilde C_\varepsilon e^{p(k+1)T},
\end{split}
\end{equation}
for a suitable $\tilde C_\varepsilon>0$. By
plugging \eqref{XY1} and \eqref{XY2} into \eqref{XY},
we complete the proof of \eqref{L3ux}
in the sub-case $\rho_k\ge\varepsilon (k+1)T$ too.

Having completed the proof of \eqref{L3ux},
we use it to complete the proof of \eqref{L3}
in the case $n=2$: indeed, by
\eqref{L3-cor-ii} and \eqref{L3ux}, when $n=2$ we have
\begin{equation}\label{775}
\ell_3 \le
{\,\rm const\,} e^{ -p(k+1)T }
\int_0^{\rho_k} 
e^{
p \sqrt{r^2+x_{o,n}^2}
} \,dr
\le {\,\rm const\,}
\big( \varepsilon^2 kT+\tilde C_\varepsilon\big).
\end{equation}
This proves \eqref{L3} also in the case $n=2$, 
since (recalling \eqref{ET}), we may take
\begin{equation}\label{Te}
T\ge 1/\varepsilon.\end{equation}
So,
the proof of \eqref{L3} is completed.

Now, we observe that
\begin{equation}\label{Sa5}
\begin{split}
\ell_1&\le{\,\rm const\,}
\int_{B_{kT}^+(x_o)} (v_k+1)^p\,dx\\ &\leq
{\,\rm const\,}
\int_{B_{kT}(x_o)} e^{p\big(|x-x_o|-(k+1)T\big)}\,dx
\\ &\leq
{\,\rm const\,}
\int_0^{{kT}} r^{n-1} e^{
p(r-(k+1)T)
}\,dr
\\ &\leq
{\,\rm const\,} (kT)^{n-1} e^{
-p(k+1)T
}
\int_0^{{kT}} e^{
pr
}\,dr\\
&=
{\,\rm const\,} (kT)^{n-1} e^{
-pT
}\\
&\leq
{\varepsilon k^{n-1}},
\end{split}
\end{equation}
provided that~$T$ is sufficiently large, possibly in
dependence of~$\varepsilon$. This last requirement, recalling
also \eqref{ET} and \eqref{Te}, fixes~$T$
once and for all (in dependence of $\varepsilon$, which, in turn,
will be fixed in
the forthcoming Lemma~\ref{H12}).

Furthermore, since~$S$ is non-decreasing and
bounded by~$1$, we obtain that
\begin{equation}\label{Sa6}
\begin{split}
&\ell_2=
\int_{ \{\theta\geq u\geq v_k
\}\cap (B_{(k+1)T}^+ (x_o)-B^+_{kT} (x_o))} 
S(v_k)\,dx
\\ &\qquad\qquad+
\int_{ \{u>\theta\}
\cap \{u\geq v_k
\}\cap (B^+_{(k+1)T}(x_o)-B^+_{kT} (x_o) )} 
S(v_k)\,dx
\\&\qquad\leq
\int_{ \{\theta\geq u>v_k
\}\cap \big(B^+_{(k+1)T}(x_o)-B^+_{kT}(x_o)\big)} 
S(u)\,dx\\ &\qquad\qquad
+{\mathcal{L}}^n\Big(
{ \{u>\theta\}
\cap \big(B^+_{(k+1)T}(x_o)-B^+_{kT}(x_o)\big)} 
\Big).
\end{split}
\end{equation}
Moreover, using also~\eqref{RES}, \eqref{61}
and~\eqref{Sa1}, and recalling \eqref{VK} once more,
we get
\begin{equation}\label{9sikjedtfyjb8882}
\begin{split}
& \int_{ B_{kT}^+(x_o)\cap\{u\le\theta\} } 
S(u)\,dx
\leq\int_{ B_{kT}^+(x_o)\cap\{ \theta\ge u>v_k\} } S(u)\,dx+
\int_{ B_{kT}^+(x_o) \cap\{ u\leq v_k\} } S(v_k)\,dx\\
&\qquad \leq{\,\rm const\,}\int_{\{\theta\ge u>v_k\}\cap\R^n_+}
F(u)\,dx+\int_{B_{kT}^+(x_o)} S(v_k)\,dx\\ 
&\qquad \leq
{\,\rm const\,}\left[
\int_{\{u>v_k\}\cap\R^n_+}
S(v_k)\,dx
+\int_{\{u>v_k\}\cap\{x_n=0\}}
S(v_k)\,dx'
+\int_{B_{kT}^+(x_o)} S(v_k)\,dx\right]
\\ &\qquad \leq 
{\,\rm const\,}(\ell_1+\ell_2+\ell_3).
\end{split}
\end{equation}
Now,
it is convenient to introduce
the following
quantities:
\begin{equation}\label{FVp12} V_r:=
\mathcal{L}^n \Big( B_r^+(x_o)
\cap\{ u>\theta \}
\Big)\qquad
{\mbox{ and }}\qquad
{A}_r:=
\int_{ B_r^+(x_o)
\cap
\{u\leq \theta \}} S(u)\,dx.\end{equation}
These quantities are appropriate variations
of similar ones defined
in \cite{CC}, and they somewhat
play the role of ``volume'' and ``area terms'', respectively,
in the minimal surface analogue.
By collecting the estimates in~\eqref{6dwrdrrq00wk1},
\eqref{9sikjedtfyjb8882},
\eqref{Sa5}, \eqref{Sa6}
and~\eqref{L3},
we conclude that
\begin{eqnarray*}
A_{kT}+
V_{kT}^{(n-1)/n}&\leq&{\rm const\,}(\ell_1+\ell_2+\ell_3)
\\ &\leq&
{\,\rm const\,} \Big[
\int_{ \{\theta\geq u\ge v_k
\}\cap (B_{(k+1)T}^+(x_o)-B_{kT}^+(x_o))} 
S(u)\,dx\\&&\quad
+{\mathcal{L}}^n\Big(
{ \{u>\theta\}
\cap (B_{(k+1)T}^+(x_o)-B_{kT}^+(x_o))} 
\Big)
+\frac{
\varepsilon k^{n-1}}{2}+C_\varepsilon (kT)^{n-2}\Big]
\\&\le&
{\,\rm const\,} \left( (
V_{(k+1)T}-V_{kT})
+(A_{(k+1)T}-
A_{kT})
+\frac{\varepsilon k^{n-1}}{2}+C_\varepsilon(kT)^{n-2}
\right).
\end{eqnarray*}
We define~$ k_\varepsilon$ to be the smallest integer
bigger than~$\mu_1+(
2C_\varepsilon T^{n-2}/\varepsilon)$, where $\mu_1$ is as in 
\eqref{eqpass0}.
This gives that $C_\varepsilon(kT)^{n-2}\le \varepsilon k^{n-1}/2$ and so
\begin{equation}\label{recurssion:21}
A_{kT}+
V_{kT}^{(n-1)/n}
\leq
{\,\rm const\,} \left( (
V_{(k+1)T}-V_{kT})
+(A_{(k+1)T}-
A_{kT})
+{\varepsilon k^{n-1}}
\right)\end{equation}
for any~$k\in \N$, with~$k\ge k_\varepsilon$.
Notice that, since~$T$ has been fixed in dependence
of $\varepsilon$ after~\eqref{Sa5},
it is conceivable to keep track of the dependence
of~$k_\varepsilon$ on~$\varepsilon$ only and disregard the
dependence on~$T$.

So, it is convenient to recall the following
general recursive result, which is a variation of
an argument in \cite{CC} and whose detailed proof may be found
in Lemma~12 of~\cite{FV}:

\begin{lemma}\label{H12}
Let $C\ge1$, $\varepsilon>0$.
Let~$\AAA_k$ and~$\VVV_k$ be two sequences
of non-negative real numbers, for~$k\in\N$.

Suppose that 
\begin{equation}\label{Bo:1}
\VVV_k \geq 1/C
\end{equation}
and
\begin{equation}\label{recurssion:2}
\VVV_{k}^{(n-1)/n}+\AAA_{k}
\leq C \Big( (
\VVV_{k+1}-\VVV_{k})
+(\AAA_{k+1}-
\AAA_{k})
+\varepsilon k^{n-1}
\Big)\end{equation}
for any~$k\in\N$.

Let
$$ c:=\min \left\{ \frac{1}{C},\,
\frac{1}{\Big( 2C (n+1)!\Big)^n }\right\}\,.
$$
Suppose that
\begin{equation}\label{FE} \e\leq \min\left\{\frac{c}{4C},\,
\frac{c^{(n-1)/n} (\sqrt[n]2-1)
}{2C}
\right\}\,.\end{equation}
Then,
\begin{equation}\label{ghdjf76738288}
\AAA_k + \VVV_k\geq c k^n
\end{equation}
for any $k\ge 4C (n+1)!$.
\end{lemma}

With this, we
define ${\mathcal{A}}_k:= A_{(k+k_\varepsilon)T}$
and ${\mathcal{V}}_k:= V_{(k+k_\varepsilon)T}$, we have that $
{\mathcal{V}}_k\ge V_{\mu_1}\ge \mu_2$, by \eqref{eqpass0},
and so \eqref{Bo:1}
holds true, if $C$ is chosen large enough.
Also, \eqref{recurssion:2}
follows from \eqref{recurssion:21}, again by choosing $C$
appropriately large.

Hence, we can exploit Lemma \ref{H12}
(notice that \eqref{FE} fixes now the value of~$\varepsilon$), 
and we deduce from \eqref{ghdjf76738288}
that
$$A_{kT}+V_{kT}
\geq\,{\rm const\,} T^n k^n$$
as long as~$k$ is large enough.

Since, by~\eqref{RES}, \eqref{61} and
\eqref{DS}, we have that 
$$ A_r\leq {\rm const\,}\int_{B_r \cap
\{u\leq \theta \} } F(u)\,dx
\leq \,{\rm const}\,r^{n-1}$$
for any~$r\ge 1$, we conclude
that~$V_r\geq{\rm const\,} r^n$
for any~$r$ suitably large,
that is \eqref{DS1}.

\bibliographystyle{amsalpha}
\bibliography{Q-bibliofile}

\end{document}